\documentclass[12pt,reqno]{amsart}
\usepackage{amsmath, amsthm, ams
symb, amsopn, amsfonts, amstext, enumerate, color, mathtools, 
mathrsfs, url}
\usepackage[shortlabels]{enumitem}
\usepackage{graphicx,color}
\usepackage{hyperref}
\usepackage[margin=1.2in]{geometry}

\makeatletter
\newcommand{\opnorm}{\@ifstar\@opnorms\@opnorm}
\newcommand{\@opnorms}[1]{%
 \left|\mkern-1.5mu\left|\mkern-1.5mu\left|
 #1
 \right|\mkern-1.5mu\right|\mkern-1.5mu\right|
}
\newcommand{\@opnorm}[2][]{%
 \mathopen{#1|\mkern-1.5mu#1|\mkern-1.5mu#1|}
 #2
 \mathclose{#1|\mkern-1.5mu#1|\mkern-1.5mu#1|}
}
\makeatother

\theoremstyle{plain}

\begin{document}


\theoremstyle{plain}
\newtheorem{theorem}{Theorem} [section]
\newtheorem{corollary}[theorem]{Corollary}
\newtheorem{lemma}[theorem]{Lemma}
\newtheorem{proposition}[theorem]{Proposition}


\theoremstyle{definition}
\newtheorem{definition}[theorem]{Definition}
\theoremstyle{remark}
\newtheorem{remark}[theorem]{Remark}

\numberwithin{theorem}{section}
\numberwithin{equation}{section}

\title[ Finite radial Morse index solutions.]{A negative result on regularity estimates of finite radial Morse index solutions to elliptic problems.}

\thanks{The authors have been supported by the FEDER-MINECO Grant PID2021-122122NBI00 and by J. Andalucia (FQM-116).}

\author[J. Silverio Mart\'{i}nez-Baena]{J. Silverio Mart\'{\i}nez-Baena}
\address{J. Silverio Mart\'{i}nez-Baena\textsuperscript{1} ---
 Departamento de An\'{a}lisis Matem\'atico, Universidad de Granada, 18071 Granada, Spain.
}
\email{jsilverio@ugr.es}

\author[Salvador Villegas]{Salvador Villegas}
\address{Salvador Villegas\textsuperscript{2} ---
 Departamento de An\'alisis Matem\'atico, Universidad de Granada, 18071 Granada, Spain.
 }
\email{svillega@ugr.es}

\date{}

\begin{abstract}
In the regularity theory of solutions to elliptic partial differential equations often the concept of stability plays the role of a sufficient condition for smoothness. It is a natural question to ask if this holds true for nonstable but finite Morse index solutions. We provide a negative answer showing the existence of sequences of solutions with radial Morse index equal to 1 for which regularity estimates can not be satisfied.
\end{abstract}

\maketitle


\section{Introduction}\label{section.intro}

In this paper we will be concern about the Morse index stability properties of solutions $u:\Omega\rightarrow \mathbb{R}$ of the nonlinear elliptic Dirichlet problem:
\begin{equation}\label{main_equation}
\left\lbrace
\begin{array}{rcll}
-\Delta u &=&  f(u)&\text{ in }\Omega\\
u&=&0&\text{ on }\partial \Omega
\end{array}\right.
\end{equation}
We will focus on $\Omega=B_1$ the unit ball in $\mathbb{R}^N$ with dimension $N\geq 3$ and the nonlinearity $f\in C^1(\mathbb{R})$ being nonnegative and nondecreasing. Under these assumptions, it is well known by the celebrated Gidas-Ni-Nirenberg symmetry result \cite{niremberg} that $u$ is radially symmetric and decreasing. For the sake of generality of the following definitions let us take $\Omega\subset\mathbb{R}^N$ to be any bounded domain and $\Omega'\subset \Omega$ a subdomain. The problem \eqref{main_equation} corresponds to the Euler-Lagrange equation for the energy functional

\begin{equation}\label{energy}
E_{\Omega}[u]=\int_{\Omega} \left(\frac{1}{2}|\nabla u|^2-F(u)\right)dx,
\end{equation}
where $F(t)=\int_0^t f(s)ds$. Consider the second variation of \eqref{energy} that when $f'(u)\in L^1_{loc}(\Omega)$ is given by  
\begin{equation*}\label{stability}
Q_{u,\Omega}[\varphi]:=\left.\frac{d^2}{d\varepsilon^2}\right|_{\varepsilon=0}E_{\Omega}[u+\varepsilon\varphi]:=\int_{\Omega} \left(|\nabla \varphi|^2-f'(u)\varphi^2\right)dx.
\end{equation*}
We say that a solution $u$ of \eqref{main_equation} is stable in $\Omega'$ if $Q_{u,\Omega'}[\varphi]\geq 0$ for all $\varphi\in C^1_0(\Omega')$ (if $\Omega'=\Omega$ we simply say that u is stable). A solution is locally stable if for any $x\in \Omega$ there exists an open neighborhood $\omega_x$ of $x$ such that $Q_{u,\omega_x}[\varphi]\geq 0$ for all $\varphi \in C_0^1(\omega_x)$. 

Quite recently, in \cite{cabre-figalli-rosoton-serra} the authors have shown that for nonlinearities $f\geq 0$ in dimensions $3\leq N\leq 9$, stable solutions to \eqref{main_equation} are as smooth as the regularity of $f$ permits, proving interior estimates:
\begin{eqnarray}
||\nabla u||_{L^{2+\gamma}(B_{1/2})} & \leq & C||u||_{L^1(B_1)}\ \ \ ( N\geq 1)\\
||u||_{C^{\alpha}(\overline{B}_{1/2})} & \leq & C||u||_{L^1(B_1)}\ \ \ (1\leq N\leq 9)\label{interiorL2gamma}
\end{eqnarray}
with $\gamma>0$, $\alpha$ and $C$ dimensional constants and analogous results up to the boundary in the case $f$ nondecreasing and convex.  
Furthermore, in \cite[Th.1.6]{villegas} the author proved sharp pointwise estimates for stable solutions of \eqref{main_equation}:
\begin{theorem}{\cite[Th.1.6]{villegas}}
Let $N\geq 2$, $f\in C^1(\mathbb{R})$, and $u\in H^1(B_1)$ be a stable radial solution of \eqref{main_equation}. Then there exists a constant $M_N$ depending only on $N$ such that:
\begin{itemize}
\item[$i)$] If $N < 10$, then $||u||_{L^\infty(B_1 )}\leq M_N ||u||_{H^1(B_1\backslash\overline{B_{1/2}})}$.\label{result_to_reproduce}
\item[$ii)$] If $N = 10$, then $|u(r)|\leq M_N ||u||_{H^1(B_1\backslash\overline{B_{1/2}})}(|\log (r)| + 1),\ \ \forall r \in (0, 1]$.
\item[$iii)$] If $N > 10$, then $|u(r)|\leq M_N ||u||_{H^1(B_1\backslash\overline{B_{1/2}})}r^{-N/2+\sqrt{N-1}+2}\ \ \ \forall r \in (0, 1]$.
\end{itemize}\label{Salvador_theorem}
\end{theorem}
These results provide a good control on the behaviour of solutions under the assumption of stability, which can be a too strong assumption in some situations so it is reasonable to require results in the same direction for a larger class of solutions. Since its first introduction by M. Morse in the early 20th century in \cite{morse}, it is now fairly well established the interest on the so called \textit{Morse index} (see for instance \cite{figalli-zhang},\cite{bahri-lions},\cite{dupaigne} and references therein), a weaker concept that contains stability as a particular case.

\begin{definition}[Morse Index $ind(u,\Omega',f)$]
Let $u$ be a solution of \eqref{main_equation} with $f'(u)\in L^1_{loc}(\Omega')$ for a subdomain $\Omega'\subset\Omega$. We say that $u$ has Morse index $ind(u)=k\in \mathbb{N}$ in $\Omega'$ if $k$ is the maximal dimension of a subspace $X_k\subset C^1_0(\Omega')$ such that,
\begin{equation*}
\int_{\Omega'}|\nabla \varphi|^2-f'(u)\varphi^2 < 0\ \ \ \ \ \ \forall \varphi\in X_k\setminus\lbrace 0 \rbrace
\end{equation*} 
If we restrict the set $C^1_0(\Omega')$ to its intersection with the space of radial functions, say $C^1_{0,rad}(\Omega')$, we get the definition of \textit{radial Morse index} in  $\Omega'$ denoted by $ind_r(u)$. As before we will omit the reference to the domain in the case $\Omega'=\Omega$. 

We shall use boundedness of radial Morse index as the criterium to the extent the class of solutions under study.
\end{definition}

\begin{remark}
It is immediate from the definition that,
\begin{itemize}
\item[$i)$] $ind(u)=0\ \text{in}\ \Omega'\iff u$ is stable in $\Omega'$.
\item[$ii)$] $ind_r(u,\Omega')\leq ind(u,\Omega')$.
\end{itemize}
\end{remark}

\begin{proposition}{\cite[Chap. 1]{dupaigne}}\label{dupaigne_proposition} Let $\Omega$ be a bounded domain and $u$ solution of \eqref{main_equation} with $f\in C^1(\mathbb{R})$. Then:
\begin{enumerate}\label{morse_index_properties}

\item \label{local_stability}
If $ind(u)< +\infty$, then $u$ is locally stable.

\item \label{finite_morse_index_for_regular_solutions}
If $u\in C^2(\overline{\Omega})$, then $ind(u)< +\infty$ and it is equal to the number of negative eigenvalues of the linearized operator $-\Delta -f'(u)$ (with Dirichlet boundary conditions).

\item $ind_r(u)=0\iff ind(u)=0$. \label{radial_morse_index_0_iff_morse_index_0}

\end{enumerate}

\end{proposition}
Proposition \ref{morse_index_properties} \eqref{radial_morse_index_0_iff_morse_index_0} is a consequence of the well known increasing behaviour of the sequence of eigenvalues of the linearized operator $(-\Delta -f'(u))$ and the fact that its first eigenvalue is simple and for $\Omega=B_1$ the corresponding eigenfuction is radial.

 Morse index had been proposed to characterize the uniform boundedness of a solution for the first time in \cite{bahri-lions}. In that paper the authors consider a class of subcritical nonlinearities, i.e., nonlinearities with asymptotic behaviour $f(t)t^{-1}|t|^{-(p-1)}\rightarrow C$ as $t\rightarrow \pm\infty$ with $1<p<\frac{N+2}{N-2}$ and prove the equivalence between boundedness of solutions and of their Morse index. In \cite{figalli-zhang} the authors achieve similar results for the supercritical case being able to prove uniform boundedness for convex and bounded domains and a suitable compact subclass of nonlinearities (positive, nondecreasing and convex). They also mention that in the critical case, one can find a counterexample to the boundedness of finite Morse index solutions. More precisely, if one define $U:\mathbb{R}^+\times\Omega\rightarrow \mathbb{R}$,
\begin{equation*}
U(\lambda,x)=\left(\frac{\sqrt{\lambda N(N-2)}}{\lambda^2+|x|^2}\right)^{\frac{N-2}{2}}
\end{equation*}
then the family of functions,
\begin{equation*}
u_\lambda(x)=U(\lambda,x)-U(1,x)
\end{equation*}
are Morse index $ind(u_\lambda)=1$ solutions (see \cite{de_marchis-pacella-ianni} Sec. 5 for a clear exposition of the details) to the problem \eqref{main_equation} for the critical nonlinearity $f_\lambda(u)=(\lambda+u)^{\frac{N+2}{N-2}}$ that exhibit a singular behaviour $\lbrace||u_\lambda||_{L^\infty(B_1)}\rbrace\rightarrow\infty$, $\lbrace||u_\lambda||_{L^1(B_1)}\rbrace\rightarrow 0$  as $\lambda\rightarrow 0^+$. This example shows that one can not expect to get something similar to Theorem \ref{result_to_reproduce} $i)$ only in terms of the finite Morse index assumption. Nevertheless, it still interesting to ask if one can control quotients like $\frac{||\cdot||_{L^p (B_1)}}{||\cdot||_{L^q(B_1)}}$ for $q<p$ under the same conditions. Our main result is again a negative answer using a bounded radial Morse index sequence as a counterexample.

\section{Main Theorem.}

\begin{theorem}\label{main_theorem}
Let $1\leq q < p \leq \infty$ with $p>\frac{N}{N-2}$. For $3\leq N\leq 9$ there exists a sequence $\{u_n\}_{n\in \mathbb{N}}\subset C^\infty(\overline{B_1})$ of solutions to Dirichlet problems \eqref{main_equation} with nonlinearities $\lbrace f_n\rbrace_{n\in \mathbb{N}}\subset C^\infty(\mathbb{R})$ and radial Morse index $ind_r(u_n)= 1\ \ \forall n\in\mathbb{N}$, such that:
\begin{equation}\label{main_theorem_equation}
\frac{||u_n||_{L^p (B_1)}}{||u_n||_{L^q(B_1)}}\longrightarrow +\infty\ \ \ \ \ \text{as}\ \ n\rightarrow\infty
\end{equation} 
\end{theorem}

\begin{remark}

It is remarkable the fact that we are not able to find divergent sequences for $q<p\leq\frac{N}{N-2}$. The ``raison d'être'' of this kind of bound in $p$ and its optimality remains as an interesting open question. We want also to emphasize again that this divergent behaviour is contrary to the stability case in which the estimate in the Theorem \ref{Salvador_theorem} $i)$ automatically implies the quotient that appears in \eqref{main_theorem_equation} to be bounded.
\end{remark}

Our main result can be interpreted as the impossibility to find a sort of reciprocal result of Prop. \ref{morse_index_properties}-\textit{(\ref{finite_morse_index_for_regular_solutions})}, i.e., the impossibility to prove certain notion of regularity of a solution only from the boundedness of its Morse index. It also remains as an open problem to find sufficient conditions to ensure such a reciprocal.
%

\section{Proof of the main theorem.}

For the sake of clarity before going into the proof of the main theorem we show the following tree lemmas.

\begin{lemma}\label{lemma1}
Let $\Omega=B_1$, $u(x)=u(|x|)\equiv u(r)$ a solution of \eqref{main_equation}. Assume $ind_r(u)=0$ in $B_{\delta}$ and in $B_{1}\setminus\overline{B_\delta}$ for some $\delta\in (0,1)$. Then $ind_r(u)\leq 1$ in $B_1$.
\end{lemma}
\noindent {\bf Proof.}
If on the contrary $ind_r(u)>1$, we could find a two dimensional subspace $X\subset C^1_{0,rad}(B_1)$ such that $Q_{u,B_1}[\varphi]<0$ for any $\varphi\in X\setminus\lbrace 0 \rbrace$. We choose two linearly independent $\varphi_1,\varphi_2\in X$. Clearly $\varphi_2(\delta)\not= 0$ because otherwise 
$$
Q_{u,B_1}[\varphi_2]=Q_{u,B_\delta}\left[\varphi_2|_{_{B_{\delta}}}\right]+Q_{u,B_1\setminus\overline{B_\delta}}\left[\varphi_2|_{_{B_1\setminus\overline{B_{\delta}}}}\right]\geq 0,
$$
which is immediately contradictory. We can define,
$$\varphi=\frac{\varphi_1(\delta)}{\varphi_2(\delta)}\varphi_2-\varphi_1$$
that vanishes at $r=\delta$. Hence
$$
Q_{u,B_1}[\varphi]=Q_{u,B_\delta}\left[\varphi|_{_{B_{\delta}}}\right]+Q_{u,B_1\setminus\overline{B_\delta}}\left[\varphi|_{_{B_1\setminus\overline{B_{\delta}}}}\right]\geq 0,
$$
which is a contradiction.\qed

\begin{lemma}{\cite[Lemma 2.1]{cabre-capella}}\label{lemma1,5}
Let $u(x)=u(|x|)=u(r)$ be a radially symmetric solution of \eqref{main_equation}. Then, $u$ is stable in $B_1\setminus\overline{B_{r_0}}$ if and only if:
\begin{equation*}
\frac{\int_{r_0}^1 r^{N-1}u_r^2\omega'^2\ dr}{\int_{r_0}^1 r^{N-1} u_r^2\frac{\omega^2}{r^2}\ dr}\geq N-1\ \ \ \ \ \ \forall \omega\in C^\infty_{0,rad}((r_0,1))\setminus\lbrace 0 \rbrace
\end{equation*}
where we have denoted $u_r:=\frac{\partial u(r)}{\partial r}$.
\end{lemma}

Actually, our Lemma \ref{lemma1,5} is an easy adaptation and not exactly the same statement of \cite[Lemma 2.1]{cabre-capella} which is slightly more general.

\begin{lemma}\label{lemma2}
Let $0<\Psi\in C^\infty ((0,1])$ such that $\Psi(t)=t, \, t\in (0, \alpha)$, for some $\alpha\in (0,1)$. Define $u:B_1\subset \mathbb{R}^N\rightarrow \mathbb{R}$, $u(x)=u(|x|)\equiv u(r)$ for $r\in [0,1]$ by
\begin{equation}
u(r)=\int_r^1\Psi(s^N)s^{1-N}ds\ \  \ \ \forall r\in[0,1].\label{definition_of_u}
\end{equation}
Then, $u\in C^\infty (\overline{B_1})$ is a solution to \eqref{main_equation} for some $f\in C^\infty(\mathbb{R})$ satisfying:

\begin{itemize}
\item[$I)$] $f\geq 0\ \text{in}\  u(B_1)\iff \Psi\label{psi_i}$ nondecreasing.\label{I)}
\item[$II)$] $f'\geq 0\ \text{in}\  u(B_1)\iff \Psi\label{psi_ii}$ is concave.
\end{itemize}
\end{lemma}
\noindent {\bf Proof.}
Since the radial function $u\in  C^\infty (\overline{B_1}\setminus\{ 0\})$ satifies $u_r<0$ for every $0<r\leq 1$ we have
$$
f(s):=-(\Delta u)({u}^{-1}(s)), \ s\in[0,u(0)) 
$$\\
is well defined and satisfies $f\in C^\infty ([0,u(0)))$. On the other hand, from $\Psi(t)=t, \, t\in (0, \alpha)$, we obtain $u_r(t)=-t$ for $t\in (0,\alpha^{1/N})$. Consequently
$$ u(r)=u(0)-\frac{r^2}{2}, \ r\in [0,\alpha^{1/N}).$$
Therefore $u\in  C^\infty (\overline{B_1})$ and $f(s)=N$ for every $s\in (u(0)-(\alpha^{1/N})^2/2,u(0))$. Finally, extending $f$ to $\mathbb{R}$ by $f(s)=N$, if $s\geq u(0)$ and in a $C^\infty$ way if $s<0$, we conclude that $f\in C^\infty(\mathbb{R})$ and $u\in C^\infty (\overline{B_1})$ is a solution to \eqref{main_equation} as claimed.
On the other hand, under radial symmetry \eqref{main_equation} is written as
\begin{equation}\label{1.1-sphericall}
-u''(r)-\frac{N-1}{r}u'(r)=f(u(r)),\, r\in (0,1].
\end{equation}
Taking the first derivative in \eqref{definition_of_u}, 
\begin{equation}\label{first_derivative_u}
-u'(r)=r^{1-N}\Psi(r^N),\, r\in (0,1].
\end{equation}
Taking the second derivative and using \eqref{1.1-sphericall} we have 
\begin{equation}\label{relation-f-psi}
f(u(r))=N\Psi'(r^N)),\, r\in (0,1],
\end{equation}
so the first part $I)$ is proven. Taking derivatives again in \eqref{relation-f-psi} and using \eqref{first_derivative_u},
$$f'(u)=-N^2r^{2(N-1)}\frac{\Psi''(r^N)}{\Psi(r^N)},\, r\in (0,1]$$
and $II)$ follows.

\

\noindent {\bf Proof of Theorem \ref{main_theorem}.}
Take an arbitrary $r_0\in (0,1)$ and define the function,
\begin{equation*}
\Psi_{r_0}(r)=\left\lbrace\begin{array}{ccc}
r         &      if       &     0< r\leq r_0^N\\
\xi_{r_0}(r)    &      if       &      r_0^N< r\leq 1
\end{array}\right.
\end{equation*}
where $\xi_{r_0}(r)$ is a $C^\infty\left([r_0,1]\right)$ strictly increasing and concave function, such that $\Psi_{r_0}\in C^\infty((0,1])$ and chosen to be bounded by 
\begin{equation}\label{bound_on_Psi}
\xi_{r_0}(r)\leq\kappa_Nr_0^N,\ r\in (r_0^N,1],
\end{equation}
with $\kappa_N=\frac{N}{2\sqrt{N-1}}$. Note that this is always possible, since $N\geq 3$ implies $\kappa_N>1$. 
Now we define a radially symmetric function $u_{r_0}:B_1\rightarrow \mathbb{R}$ as
$$
u_{r_0}(r)=\int_r^1\Psi_{r_0}(s^N)s^{1-N}\,ds\ \  \ \ \forall r\in[0,1].
$$
By the Lemma \ref{lemma2} above $u_{r_0}\in C^\infty (\overline{B_1})$ is a solution to \eqref{main_equation} for some $f_{r_0}\in C^\infty(\mathbb{R})$ satisfying $f_{r_0},f'_{r_0}\geq0$. As a result, by Proposition \ref{morse_index_properties}, $u_{r_0}$ must have finite Morse index and be locally stable. Indeed $u_{r_0}$ has zero radial Morse index in $B_{r_0}$ since by \eqref{psi_ii} and \eqref{relation-f-psi} $f_{r_0}(u_{r_0}(r))=N$ and thus $f_{r_0}'(u_{r_0}(r))=0$ for $r<r_0$. It follows immediately that  $u_{r_0}$ is stable in $B_{r_0}$. Taking into account Lemma \ref{lemma1,5} we can prove stability also in the annulus $B_1\setminus\overline{B_{r_0}}$:
\begin{eqnarray*}
\frac{\int_{r_0}^1 r^{N-1}{(u_{r_0})_r}^2\omega'^2\ dr}{\int_{r_0}^1 r^{N-1}{(u_{r_0})_r}^2\ \frac{\omega^2}{r^2}\ dr}&=&
 \frac{\int_{r_0}^1 r^{N-1}(r^{2-2N}\Psi_{r_0}^2(r^N))\omega'^2\ dr}{\int_{r_0}^1 r^{N-1}(r^{2-2N}\Psi_{r_0}^2(r^N))\frac{\omega^2}{r^2}\ dr}\geq  \frac{\int_{r_0}^1 r^{N-1}(r^{2-2N}r_0^{2N})\omega'^2\ dr}{\int_{r_0}^1 r^{N-1}(r^{2-2N}\kappa_N^2 r_0^{2N})\frac{\omega^2}{r^2}\ dr}\geq\\
& \geq & \frac{r_0^{2N}}{\left(\frac{N}{2\sqrt{N-1}}r_0^N\right)^2}\frac{N^2}{4}
 =  N-1
\end{eqnarray*}
where we have used \eqref{bound_on_Psi} and the generalized Hardy inequality ($0<a<b$),
 $$\int_a^b r^{\alpha+1}\omega'^2\ dr\geq\frac{\alpha^2}{4}\int_a^b r^{\alpha-1}\omega^2\ dr\ \ \ \forall \alpha\in \mathbb{R},\ \forall \omega\in C^1_0([a,b])$$
  had been used with $a=r_0,b=1$ and $\alpha=-N$. Then Lemma \ref{lemma1} implies that the radial Morse index $ind_r(u_{r_0})\leq 1$. 
  
  Now, we compute the estimates for its $L^p$ and $L^q$ norms in the case $\frac{N}{N-2}<p<\infty$. Note that, in this case, there is no loss of generality in assuming $N/(N-2)<q<p$ since $L^{q_2} (B_1)$ is continuously embedded in $L^{q_1} (B_1)$ if $q_1<q_2$.
\begin{eqnarray}\nonumber
||u_{r_0}||^p_{L^p(B_1)}&=&C_N\int_0^1r^{N-1}\left|\int_r^1s^{1-N}\Psi_{r_0}(s^N)ds\right|^p
dr\\\nonumber
&\geq &C_N\int_0^{r_0} r^{N-1}\left(\int_r^{r_0}s\,ds\right)^p dr= C_N\int_0^{r_0} r^{N-1}\left( \frac{r_0^2-r^2}{2}\right)^p dr \\\nonumber
&= &C_{N,p}\ r_0^{N+2p}.\\\nonumber\\\nonumber
||u_{r_0}||^q_{L^q(B_1)}
&= & C_N\int_0^{r_0}r^{N-1}\left(		\frac{r_0^2-r^2}{2}+\int_{r_0}^1 s^{1-N}\Psi_{r_0}(s^N)ds\right)^q\,dr\\\nonumber
&+&C_N\int_{r_0}^1r^{N-1}\left(\int_r^1s^{1-N}\Psi_{r_0}(s^N)ds\right)^q dr\\\nonumber
&\leq& C_N\int_0^{r_0}r^{N-1}\left(		\frac{r_0^2-r^2}{2}+\kappa_n r_0^N\int_{r_0}^1 s^{1-N}ds\right)^q\,dr\\\nonumber
&+&C_N\int_{r_0}^1r^{N-1}\left(\kappa_n r_0^N \int_r^1s^{1-N}ds\right)^q dr\\\nonumber
&\leq& C_N\int_0^{r_0}r^{N-1}\left(		\frac{r_0^2-r^2}{2}+\kappa_n r_0^N \frac{r_0^{2-N}}{N-2}\right)^q\,dr\\\nonumber
&+&C_N\int_{r_0}^1r^{N-1}\left(\kappa_n r_0^N \frac{r^{2-N}}{N-2}\right)^q dr\\\nonumber
&\leq&C_{N,q}r_0^{2q}\int_0^{r_0}r^{N-1}dr+C_{N,q}r_0^{Nq}\frac{r_0^{N+(2-N)q}}{-N+(N-2)q}=C'_{N,q}r_0^{N+2q}. \\\nonumber 
\end{eqnarray}
Combining these inequalities we obtain
\begin{equation}
\frac{||u_{r_0}||_{L^q(B_1)}}{||u_{r_0}||_{L^p(B_1)}}\leq C_{N,p,q}r_0^{N\left(\frac{1}{q}-\frac{1}{p}\right)}\label{quotient_ending_proof}
\end{equation}
with $C_{N,p,q}$ a dimensional constant depending on $q$ and $p$. Since $q<p$ the quotient \eqref{quotient_ending_proof} goes to $0$ as $r_0\rightarrow 0$. This proves \eqref{main_theorem_equation} in the range $\frac{N}{N-2}<p<\infty$. If $p=\infty$ we see at once that
\begin{eqnarray*}\nonumber
||u||_{L^\infty (B_1)}=u(0)=\int_0^1 s^{1-N}\Psi(s^N)\,ds\geq \int_0^{r_0} s^{1-N}\Psi(s^N)\,ds=\frac{r_0^2}{2}\\\nonumber
\end{eqnarray*}
Now we get,
\begin{equation}\label{Linfty_Lq_estimate}
\frac{||u_{r_0}||_{L^\infty(B_1)}}{||u_{r_0}||_{L^q(B_1)}}\geq C_{N,q}\,r_0^{-\frac{N}{q}}
\end{equation}
 which is unbounded as $r_0\rightarrow 0$.
 
 Finally, by Proposition \ref{dupaigne_proposition} \eqref{radial_morse_index_0_iff_morse_index_0}, if $ind_r(u_{r_0})=0$ in $B_1$ then $u$ would be stable and, taking $q=1$,\eqref{interiorL2gamma} would imply the quotient $\Vert u_{r_0} \Vert_{\infty} / \Vert u_{r_0}\Vert_1$ to be uniformly bounded by a dimensional constant, fact that we have just proven impossible for small $r_0>0$. This, together with $ind_r(u_{r_0})\leq 1$ deduced above, shows that $ind_r(u_{r_0})=1$ for sufficiently small $r_0>0$ and the theorem follows.
\qed

 \medskip

\noindent{\bf Acknowledgements:} The first author wishes to express his thanks to Xavier Cabr\'e  for many stimulating conversations during the preparation of the paper.


%


\



\end{document}